\documentclass[12pt]{amsart}
\usepackage{amsmath}
\usepackage{amsxtra}
\usepackage{amscd}
\usepackage{amsthm}
\usepackage{amsfonts}
\usepackage{amssymb}
\usepackage{eucal}
\usepackage{xypic}

\advance\textheight by -1\footskip \relax

\usepackage{verbatim} 
\usepackage{mathrsfs}

\newtheorem{cor}[subsection]{Corollary}
\newtheorem{lem}[subsection]{Lemma}

\newtheorem{thm}[subsection]{Theorem}
\newtheorem*{thmno}{Theorem}
\newtheorem*{thm*}{Theorem}
\newtheorem{defn}[subsection]{Definition}
\theoremstyle{definition}
\newtheorem{exm}[subsection]{Example}
\newtheorem{rem}[subsection]{Remark}

\newtheorem*{defn*}{Definition}

\newcommand{\thmref}[1]{Theorem~\ref{#1}}
\newcommand{\secref}[1]{Section~\ref{#1}}
\newcommand{\lemref}[1]{Lemma~\ref{#1}}

\newcommand{\defnref}[1]{Definition~\ref{#1}}

\newcommand{\eqnref}[1]{{(\ref{#1})}}

\newcommand{\nc}{\newcommand}

\renewcommand{\AA}{{\mathbf{A}}}
\newcommand{\CC}{{\mathbf{C}}}
\newcommand{\ZZ}{{\mathbf{Z}}}

\newcommand{\gf}{{\mathfrak{g}}}

\newcommand{{\gl}}{{\mathfrak{gl}}}
\newcommand{{\spl}}{{\mathfrak{sl}}}

\newcommand{\ra}{{\longrightarrow}}

\newcommand{\PP}{{\mathbf P}}

\newcommand{\zf}{{\mathfrak z}}
\newcommand{\tf}{{\mathfrak t}}

\newcommand{\bof}{{\mathfrak b}}
\newcommand{\uf}{{\mathfrak u}}

\renewcommand{\O}{{\mathcal O}}
\newcommand{\isom}{\cong}

\DeclareMathOperator{\Isom}{Isom}
\DeclareMathOperator{\Hom}{Hom}

\DeclareMathOperator{\Ad}{Ad}
\DeclareMathOperator{\Aut}{Aut}
\DeclareMathOperator{\Id}{Id}
\DeclareMathOperator{\rank}{rank}
\DeclareMathOperator{\Bun}{Bun}
\DeclareMathOperator{\Pic}{Pic}

\DeclareMathOperator{\codim}{codim}
\DeclareMathOperator{\coker}{coker}
\renewcommand{\div}{{\textnormal{div}}}

\renewcommand{\H}[1]{{\textnormal H}^{#1}}
\newcommand{\h}[1]{{\textnormal h}^{#1}}

\newcommand{\XX}{{\mathbf X}}
\newcommand{\YY}{{\mathbf Y}}

\newcommand{\lat}{{\mathbf\Lambda}} 
\newcommand{\cgr}{{\mathbf G}} 
\newcommand{\cgrr}{{\cgr_r}}

\newcommand{\cborel}{{\mathbf B}}
\newcommand{\ccartan}{{\mathbf T}}
\newcommand{\cuni}{{\mathbf U}}
\renewcommand{\a}{{\mathbf{\alpha}}}

\newcommand{\e}{{\mathbf e}}
\newcommand{\rlat}{{\mathbf Q}}
\newcommand{\lcan}{{\mathbf{\kappa}}} 
\newcommand{\rsys}{{\mathbf{\Delta}}} 
\newcommand{\rsim}{{\mathbf{\Sigma}}} 

\newcommand{\E}{{\mathbf E}}

\newcommand{\higgsp}[1]{{\textnormal{Higgs}'_{ \widetilde{#1} } (#1)} }

\newcommand{\antiw}[1]{{{\omega^{-1}_{#1}}}}

\newcommand{\W}[1]{W(#1)} 
\newcommand{\Ws}[1]{{W_e(#1)}} 

\newcommand{\wX}{{\widetilde X}}

\newcommand{\R} {{\Delta}}    
\newcommand{\re}{{\R_e}} 
\newcommand{\ri}{{\R_i}} 

\newcommand{\on}{\operatorname}
\newcommand{\Tors}{{\on{Tors}}}
\newcommand{\Q}{{\mathcal Q}}

\newcommand{\Tot}{\on{Tot}}
\newcommand{\Eig}{{\mathcal{E}}}
\newcommand{\cf}{{\mathfrak c}}

\newcommand{\subc}{{\cf}} 

\newcommand{\wlat}{{\mathbf P}}
\newcommand{\Lie}{\on{Lie}}
\newcommand{\lra}{{\longrightarrow}}

\nc{\wXX}{{\widetilde\XX}}
\nc{\wYY}{{\widetilde\YY}}
\nc{\wB}{{\widetilde B}}
\renewcommand{\SS}{{\mathbf S}}

\newcommand{\DD}{{\mathbf D}}
\newcommand{\wDD}{{\widetilde\DD}}
\newcommand{\UU}{{\mathbf U}}

\newcommand{\creg}{{\mathbf R}}

\newcommand{\keq}{{=}} 
\newcommand{\extn}[3]{{ {#1}_{#2} \times^{#2} {#3}}}
\newcommand{\cl}{{\on{cl}}}

\newcommand{\Gr}{{\on{Gr}}}
\newcommand{\rf}{{\mathfrak r}}
\newcommand{\GNb}{{ \overline{G/N} } }
\newcommand{\GTb}{{ \overline{G/T} } }
\newcommand{\Sch}{{\on{Sch}}}
\newcommand{\A}{{\mathcal A}}
\newcommand{\wY}{{\widetilde Y}}
\newcommand{\sra}{{\rightarrow}}

\newcommand{\Higgs}{\on{Higgs}}

\begin{document}
\title[Regular Bundles]{Almost Regular Bundles on del Pezzo Fibrations}
\author{K\"ur\c{s}at Aker}
\address{
Universit\"{a}t Hannover, Institut f\"{u}r Mathematik,
Lehrgebiet C, 
D-30167 Hannover, Germany 
}
\email{aker@math.uni-hannover.de}
\thanks{Supported by NSF grants DMS-01004354, FRG-0139799 and Schwerpunktprogramm ''Globale Methoden in der komplexen Geometrie''
HU 337/5-2.
}
\begin{abstract}
This paper is devoted to the study of a certain class of principal bundles on del Pezzo surfaces, which were introduced and studied by  
Friedman and Morgan in \cite{FMdP}: The two authors showed that there exists a
unique principal bundle (up to isomorphism) on a given (Gorenstein) del Pezzo surface
satisfying certain properties. 
We call these bundles
{\em almost regular}. In turn, we study them in families.
In this case, the existence and the moduli
of these bundles are governed by the cohomology groups of an abelian sheaf ${\mathscr A}$: On a given del Pezzo fibration, the existence of an almost regular
bundle depends on the vanishing of an
obstruction class in $\H{2}({\mathscr A})$. In which case,
the set of isomorphism classes of almost regular bundles 
become a homogeneous space under 
the $\H{1}({\mathscr A})$ action.
\end{abstract}
\maketitle



A del Pezzo surface, $S$, is a smooth complex projective surface whose anticanonical bundle $\antiw{S}$ is nef and big. The classification of del Pezzo surfaces shows that such a surface is either
isomorphic to $\PP^1\times \PP^1$ or is the blow-up of $\PP^2$ at $0\leq r \leq 8$ points in almost 
general position (position presque-g\'en\'erale \cite{Dem}). 
A Gorenstein del Pezzo surface $Y$ is a normal rational projective surface  whose anticanonical sheaf  $\antiw{Y}$ is invertible (=Gorenstein) and  ample. Such a surface
is the anticanonical model of a del Pezzo surface $S$. 
A principal bundle is a generalization of a vector bundle in which the fibers of the bundle, previously copies of a fixed vector space $V$, are now replaced with the copies of a fixed (complex Lie) group $G$. In this paper, we study a class of principal bundles 
on families of del Pezzo surfaces which are {\em ``natural''} in a certain sense. We will call these bundles 
{\em almost regular}.

Friedman and Morgan \cite{FMdP} construct such bundles on a single surface $Y$ and show that
they are all isomorphic. In other words, there exists
a tautological isomorphism class of bundles on $Y$. Our question 
is to study these bundles on a family of (Gorenstein) del Pezzo surfaces, $p:\YY \ra \XX$, where $p$ is a projective integral flat map whose geometric fibers are (Gorenstein) del Pezzo surfaces. We formulate our solution in terms of  
tools
developed for principal bundles on elliptic fibrations and Higgs bundles, such as cameral covers and
abstract Higgs bundles, some basic singularity theory (simultaneous resolutions of rational double points), and
(sub)regular elements from reductive Lie groups.



Our solution, most satisfactory when the base of the family is a curve (\thmref{thm:MainCurve}), shows that 
almost regular bundles on $\YY$ can be classified mainly in terms of maximal torus bundles
on a cameral cover $\wXX$ (\defnref{defn:CamCover}) of $\XX$ (\thmref{thm:MainGeneral}). Such a classification
is along the classical lines of spectral covers, where Higgs bundles on a given
variety are classified in terms of line bundles on spectral covers of the original variety (briefly reviewed in \secref{ssec:SpectralCovers}). 

The classification of almost regular bundles up to isomorphism goes as follows:
First of all, they are locally isomorphic with respect 
to the base $\XX$ (\thmref{thm:IsomLocalModel}). In other words, given two such 
bundles, for any point $x\in \XX$, there is an etale neighborhood $\UU\ra \XX$ of $x \in X$ s.t. on $
\UU \times_\XX \YY$, the bundles are isomorphic. Combined with the fact
that the automorphism groups of almost regular bundles along fibers are abelian
shows that there exists an Abelian sheaf $\mathscr{A}$ of such automorphisms:
One sees that the existence of almost regular bundles on $\YY$ are contolled by
an obstruction class in $ \H{2}(\XX,\mathscr{A})$. Once this class vanishes,
the isomorphism classes of almost regular 
bundles become a homogeneous space under the
action of $\H{1}(\XX,\mathscr{A})$. This agrees with the results of \cite{DG} and \cite{FM1}.




This paper adds yet another link to the chain of results relating del Pezzo surfaces and group theory.
Here are just two of many examples for  the close ties between del Pezzo surfaces and group theory: 
\begin{itemize}
\item
Given a smooth cubic $S$, its second cohomology $\H{2}(S,\ZZ)$ contains a copy of the root lattice of $E_6$. 
\item
The twenty seven lines of the cubic correspond to the weights of the fundamental 
representation of $E_6$.  
\end{itemize}

Observations like these as well as an impetus coming from mathematical physics (F-theory vs. heterotic string) is behind the desire
to study the questions  we do in this paper.

\noindent {\bf Conventions}. In this paper, by a scheme, we mean deal with a scheme.
Unless otherwise mentioned, we will utilize the etale topology on such a scheme. Our results are valid for complex topology as well. 

\noindent {\bf Acknowlegments} 
This paper generalizes earlier results of my Ph.D. thesis drafted under the supervision of Ron Donagi at the University of Pennsylvania. 
I would like to thank my advisor, Ron Donagi, for
suggesting this problem to me as well as 
 his constant support during my graduate studies. Tony Pantev and Sukhendu Mehrotra have contributed to this work with many stimulating conversations. I am also most grateful to Klaus Hulek for his hospitality which I have enjoyed during my stay at the University of Hannover, whose suggestions also improved the reading of the current manuscript.

\section{Del Pezzo Surfaces} \label{ssec: Del Pezzo}
Here we review some classical \cite{Dem}, \cite{Dol}, \cite{Manin} and less classical
\cite{FMdP} facts about del Pezzo surfaces.

\begin{defn}{\cite{Dol}}
A del Pezzo surface $S$ is a smooth complex projective surface whose anticanonical bundle $\antiw{S}$ is nef and big.
\end{defn}

A line bundle $L$ on a smooth surface $S$ is nef if $L\cdot C\geq 0$ for all irreducible
curves $C\subset S$. A nef line bundle $L$ on $S$ is {\em big} if $L\cdot L>0$.

\begin{defn}
A normal projective surface $Y$ is called a Gorenstein surface if the dualizing sheaf $\omega_Y$ on $Y$ is invertible. 
\end{defn}

A Gorenstein surface $Y$ with ample anticanonical bundle $\antiw{Y}$ is either rational or a cone over
an elliptic curve. The singularities on such a surface are either rational double points (in the case $Y$ is rational),
or the unique simple elliptic singularity (in the case $Y$ is a cone over elliptic curve) \cite{HW}.

\begin{defn}
A normal projective surface $Y$ is called a Gorenstein del Pezzo surface if it is a rational Gorenstein surface with $\antiw{Y}$ ample.
\end{defn}

A del Pezzo surface is either
isomorphic to $\PP^1\times \PP^1$ or is the blow-up of $\PP^2$ at $0\leq r \leq 8$ points in almost 
general position  (position presque-g\'en\'erale \cite{Dem}).
The minimal resolution of a Gorenstein del Pezzo surface is a del Pezzo surface. In return, the anticanonical model of a del Pezzo surface is a Gorenstein del Pezzo surface. For a given del Pezzo surface $S$, the anticanonical model of $S$ is defined to be $Y:=\textnormal{Proj}(\oplus_{n}
\Gamma(S, \omega_S^{-\otimes n} ) )$. 
We will call such a pair $(S,Y)$  a {\em del Pezzo pair}.


The degree of a del Pezzo surface $S$ is the number $K_S^2$ and will be denoted by $d$. The degree of a Gorenstein del Pezzo is defined accordingly. When $S\not\cong \PP^1\times\PP^1$,
then $d=9-r$. In the rest of this work, $\PP^1\times\PP^1$ will be omitted as a  (Gorenstein) del Pezzo surface.

A del Pezzo surface $S$ is simply connected. Its  second integer cohomology $\H{2}(S,\ZZ)$ carries  most of the topological information of $S$ and is also isomorphic to $\Pic(S)$. It carries a unimodular bilinear form, namely {\em the intersection pairing}. 
The following lattice is an abstract model for $\H{2}(S,\ZZ)$ for a surface of degree $9-r$:

\subsection{Unimodular lattice $\lat_r$} \label{ssec:Lattice}
Define $\lat_r=\ZZ^{1+r}$ with the (orthonormal) basis $\e_0,\ldots,\e_r$ 
satisfying $\e_i\cdot\e_j=0$ if $i\neq j$, $\e_0^2=1$ and $\e_i^2=-1$ for
$i=1,\ldots,r$ for $1\leq r\leq 8$. Inside $\lat_r$, fix
$\lcan_r:=-3\e_0+\sum_{i=1}^r\e_i$. For convenience, define $(\e_i,\e_j):=-\e_i\cdot \e_j$.
Since $\lat_r$ is unimodular, we identify $\lat_r$ by its dual via
the isomorphism induced by the bilinear form $(\cdot,\,\cdot)$.

The second cohomology $\H{2}(S,\ZZ)$ of a surface $S$ of degree $9-r$ 
is isomorphic to $\lat_r$. Among all such isomorphisms, we single out the ones with 
geometric significance:

\begin{defn}
An isometry $\varphi : \lat_r \ra \Pic(S)  $
with $\varphi(\lcan_r)=K_S$ is called a {\em marking}.
If this marking is induced by a blow-down structure on $S$, it is called a {\em 
geometric} marking.
\end{defn}

\begin{defn}[\cite{Dol}] \label{defn:BlowDownStr}
A {\em blowing down structure} on a Del Pezzo $S$ of degree $9-r$ 
is a composition
of birational morphisms 
$$\xymatrix{S=S_r\ar[r]^{\pi_r} & S_{r-1}\ar[r]^{\pi_{r-1}} & \cdots
\ar[r]^{\pi_{2}} & S_1 \ar[r]^{\pi_{1}} & S_0:=\PP^2},$$ 
where each $\pi_i: S_i \lra S_{i-1}$ is the blow-up of a point $x_i\in
S_{i-1}$ for $i=1,\ldots,r$.
\end{defn}

A blowing down structure on $S$ induces a basis of $\Pic(S)$ formed by 
the classes $e_0=\pi^*(\ell)$, $e_i=(\pi_k\circ\cdots\pi_i)^*(E_i)$,
where $\ell$ is a line on $\PP^2$ and $E_i=\pi_i^{-1}(x_i)$. 
Such a basis of $\Pic(S)$ is called 
{\em geometric basis}. The isometry sending $\e_i\in\lat_r \mapsto e_i\in\Pic(S)$ is the induced marking.

Earlier in the introduction, we have mentioned the close ties between del Pezzo surfaces and group theory. We formulate few, equivalently, in terms of $\lat_r$:
First, extend the series of simply connected group $E_6,E_7,E_8$ for $r=6,7,8$
via setting $E_r=A_1\times A_2, A_4, D_5$ for $r=3,4,5$ and denote the root and the weight lattices of $E_r$ by $\rlat(E_r)$ and $\wlat(E_r)$ respectively. Then:
\begin{itemize}

\item
The dual lattice pair $(\ZZ\lcan_r)^\perp$ and
$\lat_r/(\ZZ\lcan_r)$ is isomorphic to the lattice pair $\rlat(E_r)$ and $\wlat(E_r)$. Denote the former pair by $\rlat(\lat_r)$ and $\wlat(\lat_r)$.

\item Fix a maximal torus $T\subset E_r$ and let $W:=N_{E_r}(T)/T$ be the corresponding Weyl group. Then the group of self-isometries of $\lat_r$ preserving
$\lcan_r$ is isomorphic to $W$. The group $W$ acts simply 
transitively on the set of markings of a surface S of degree $9-r$.

\end{itemize}

Roots in $\lat_r$ are characterized as follows: $\alpha\in \lat_r$ is a root if 
$\alpha^2=-2$ and $\alpha\cdot \lcan=0$. A root $\alpha$ is {\em effective} with 
respect to a marking $\varphi$ if the corresponding class is an effective divisor on $S$. One can accordingly define (effective) roots on $S$. Irreducible effective roots on $S$ are also called $(-2)$-curves, i.e. irreducible curves $C\subset S$ s.t. $C^2=-2$ and $C\cdot K_S=0$. For a given del Pezzo pair $(S,Y)$, the birational 
morphism $S\ra Y$ contracts
exactly the $(-2)$-curves. The resulting singularities are rational double points. 
Reading this map in the reverse order, we see that $S\ra Y$ is the blowup of 
rational double points of $Y$ whose exceptional divisors are unions of special configurations
of $(-2)$-curves.
Denote the set of $(-2)$-curves, effective roots and roots in $\Pic(S)$
respectively by $\re(S)\subset \ri(S) \subset \R(S)$.

For a del Pezzo pair $(S,Y)$,
we can relate the rational double points of  $Y$,  $S$ and root systems as follows: 
Inside the root system $\R(S)$, the set of
effective curves $\ri(S)$ form a subroot system. The set of $(-2)$-curves $\re(S)$ is a set of simple roots for this subroot system. 

\subsection{Rational Double Points} \label{ssec:RDP} Rational double points are
the simplest of all isolated surface singularities:
\begin{defn}
A point $y\in Y$ is a {\em rational double point(=RDP)} if for one (hence
for all) resolution(s) $\rho:S \ra Y$, $R^1\rho_*\O_S=0$ near $y$.
\end{defn}
Like del Pezzo surfaces, rational double points on a normal surface
are also related to Lie theory: The exceptional divisor of the minimal resolution of a rational double point singularity is a connected union of  
$(-2)$-curves given by a special configuration.
This configuration is best explained via {\em the dual graph} of the singularity:
 Whenever two such curves meet, the intersection is transversal. The dual graph of the singularity is the undirected graph whose vertices are $(-2)$-curves and edges link the vertices  corresponding to intersecting $(-2)$-curves. The analytic isomorphism
class of a rational double point is determined by its dual graph, which in turn is a Dynkin diagram of type A, D, or E. For this reason, these singularities are also called ADE singularities.
Furthermore, one can construct a rational double point starting from the Lie algebra of the same type. Indeed, one can construct more, including the semi-universal deformation space and the simultaneous resolution of the rational double point singularity. These singularities have many different characterizations \cite{Durfee}. We refer the reader to \cite{Slodowy1}, \cite{Slodowy2}, \cite{Slodowy3} for different aspects of these singularities. 

\subsection{Simultaneous Resolutions}
Rational double points are exceptional among other surface singularities because
they admit simultaneous resolutions, i.e. given a deformation $Y$ of (a germ of) a rational double point
singularity over some base $X$ (again a germ), there exists a finite base change $X'\sra X$ and  a morphism $S\sra Y':=Y\times_X X'$   s.t. the diagram 
$$\xymatrix{
  S\ar[r]^r \ar[d]_s 
&  Y'\ar[d]_{\tilde{p}} \ar[r]^{\widetilde c} & Y \ar[d]_p \\
 X' \ar@{=}[r] & X' \ar[r]^c & X \\ 
},$$
commutes, 
$S \sra X'$ is smooth and $S \sra Y'$ is proper.  Here a given {\em fiber} of $S$ is a resolution of singularities of the corresponding {\em fiber} of $Y'$. 


The existence of such a local resolution was proved by Brieskorn \cite{Bri1}, \cite{Bri2} and Tjurina \cite{Tjurina1}, \cite{Tjurina2}. Later, Artin gave a global
version of this theorem \cite{Artin}: Namely, given a flat family of surfaces $Y$ over
a scheme $X$ whose fibers contain at worst rational double points, there exists a finite base change $X'\sra X$ and  $S\sra Y':=Y\times_X X'$   s.t. $S \sra X'$ is smooth and $S \sra Y'$ is proper (\cite{KM}). 

Briefly, we touch on the group theoretic nature of rational double points singularities.
W.l.o.g. assume that we only have a unique singular point. Let $\gf$ be the corresponding, simple Lie algebra, $\tf$ a Cartan subalgebra in $\gf$ and $W$ its Weyl group. Then one can choose a
$2+\rank \gf$ dimensional plane $Y$ in $\gf$ s.t. after setting $X'=\tf $ and $X=\tf/W$, we see that
$Y \sra X$ is the semiuniversal deformation space, where the fiber above $0\in X$
is the singularity of the required type. 
Furthermore,  after the base change $X' \sra X$, one can resolve the singularities of $Y'$ simultaneously.  

Unfortunately, the group theoretic 
character of the local construction is lost in the global version. Namely, in the local case, the cover $X' \sra X$ could be chosen Galois under $W$ as we did above,
 whereas in the global version there is no 
such guarantee. We will single out such simultaneous resolutions:


For a Cartan subalgebra $\tf$ (of some reductive Lie algebra $\gf$) and its Weyl group $W$, we introduce {\em $(W-)$cameral covers}. The significance of such
covers is that they can be used to construct principal bundles, more specifically, abstract Higgs bundles (\defnref{defn:AbsHiggs1}). This point will become
clearer in \secref{chap:Bundle}.

\begin{defn} \label{defn:CamCover}
A finite flat cover $\wX \ra X$ is called a ($W$-)cameral
cover if locally (in the etale topology) it is a pullback of the covering $\tf \ra \tf/W$.
\end{defn}

\begin{defn} For a family of Gorenstein del Pezzo surfaces of degree $9-r$, $\YY \ra \XX$, 
a {\em cameral resolution} of $p:\YY \ra \XX$
is a family of del Pezzo
surfaces $s:\SS\ra\wXX$, together with a map $r:\SS\ra\YY$ so that the diagram
$$\xymatrix{
  \SS \ar[r]^r \ar[d]_s 
&  \wYY\ar[d]_{\tilde{p}} \ar[r]^{\widetilde c} & \YY \ar[d]_p \\
 \wXX \ar@{=}[r] & \wXX \ar[r]^c & \XX\\ 
},$$
commutes and 
$c:\wXX \ra \XX$ is a $W$-cameral cover and $r: \SS \ra \wYY$ is a 
simultaneous
resolution over $\wXX$.
\end{defn}

Here, $W$ is the group of self isometries of $\lat_r$ preserving $\lcan_r$ (or, equivalently, the Weyl group of a Cartan subgroup of $\cgrr$ (see \secref{ssec:ConformalGroups})).

It is clear that such resolutions will not exist in general, however,  Brieskorn's earlier work \cite{Bri1}, \cite{Bri2} combined with \cite{Slodowy2}
shows that when
$\XX$ 
is a smooth projective curve and $\YY$ a smooth projective threefold, 

\begin{thm} \label{thm:CurveCamCover}
There is a cameral resolution of $\YY \ra \XX$ so that
\begin{itemize}
\item the cameral cover $\wXX$ is smooth,
\item the pullback family $\wYY$ is normal and 
\item the family of del Pezzo surfaces, $\SS$, is a smooth threefold.
\end{itemize}
\end{thm}

\subsection{Conformal Groups} \label{ssec:ConformalGroups}
Given a del Pezzo surface of degree $S$ of degree $9-r$, the sublattice $(\ZZ K_S)^\perp$ is isomorphic to the root lattice of simply connected group $E_r$, in other words to the cocharacter lattice of $E_r$, for $3\leq r\leq 8$. 
Is there another group related to $E_r$ whose cocharacter lattice is isomorphic to 
the full lattice $\H{2}(S,\ZZ)$ ? The answer is yes. These groups are constructed in
\cite{FMdP} and are called {\em conformal groups}. For a given $r$, a conformal group is of the form 
$E_r \times_{\ZZ/(9-r)\ZZ} \CC^*$ and it is essentially unique up to isomorphism \cite{FMdP}.
Denote the resulting group by $\cgrr$. We will drop $r$ if it is clear from
the context.

Conformal groups have been previously introduced to construct  universal {\em regular} bundles on the moduli space of 
principal $G$-bundles on a smooth elliptic curve for a given simple, simply connected group $G$ \cite{FM-Auto} (for regular bundles, see \defnref{defn:RegularElliptic}).


For our purposes, it will be sufficient to know the root theoretic properties, a set of fundamental weights and the following facts about the group $\cgrr$: The cocharacter lattice (of a maximal torus) of a conformal group $\cgrr$  is  isomorphic to $\lat_r$. This lattice is unimodular, thus it is isomorphic
to the character lattice. Thus we identify these two lattices. Under this identification, the roots are identified with coroots. Hence, $\cgrr$ is Langlands self-dual. In addition, the derived subgroup of $\cgrr$ 
is $E_r$. In particular, it is a simply connected semisimple group. Using these information, we will prove an important property of $\cgrr$: The centralizer of 
a given regular element $g$ in $\cgrr$ is connected (\thmref{thm:RegCentIsConnected}).




\subsection{Regularity} Regular and subregular elements of Lie groups (resp. algebras)
and related concepts are the common core of several notions in this paper.  For a given reductive group,

\begin{defn} \label{defn:Regular}
An element $g$ of a reductive Lie group is said to be {\em regular}(resp. {\em subregular}) if its centralizer
is of dimension $\rank G$ (resp. $2+\rank G$). 
\end{defn}

An element of $ GL(n)$ (resp. $\mathfrak{g}\mathfrak{l}_n$) is regular iff there
exists a unique Jordan block for each eigenvalue iff the minimal polynomial is equal to
the characteristic polynomial.

One can think of the centralizer $C_G(g)$ as the automorphism group of $g$  and $\cf_\gf(g)$ as its Lie algebra of infinitesimal endomorpshims  under
the conjugation action on $G$ .

The analogous definition on a smooth elliptic curve is
\begin{defn}
A semistable principal $G$-bundle $E_G$ on a smooth elliptic curve is {\em regular}
if $h^0(Ad(E_G))=\rank G$ (or equivalently, $h^0(Aut(E_G))=\rank G$).
\end{defn}

Indeed, \cite{FM1} links together these two concepts showing that the (infinitesimal) automorphism group
of a semistable principal $G$-bundle on a smooth elliptic curve is isomorphic to
the (infinitesimal) centralizer of a Lie algebra element. (Sub)regular elements
will also appear in the context of rational double points.

The centralizer of a regular element will be called a regular centralizer. A regular
section of a Lie algebra bundle or a group scheme is a section whose values are regular
in the above sense. 

We now want to prove the following fact about conformal groups:

\begin{thm} \label{thm:RegCentIsConnected}
All regular centralizers in $\cgrr$ are connected. 
\end{thm}

Notice that the same statement is true for a semisimple group $G$ if and only if
it is both simply connected and adjoint. As a corollary of the theorem, we see that
a regular centralizer in $\cgrr$ is a connected Abelian group of dimension $1+r$.

 Before proving this assertion, we will review some facts and results adapted
from \cite{Humphreys}:

Given an element $x$ in a Lie group $G$, denote its Jordan decomposition 
$x=x_s x_u$ where  $x_s$ is the semisimple and $x_u$ is the unipotent part of $x$.
 Two Lie group elements $x$ and $y$ commute  if 
and only if all parts of the Jordan decompositions of $x$ and $y$ commute. Then,

 \begin{lem}[Centralizer Reduction] \label{lem:CentReduction}
$$ C_G(x)=C_{C_G(x_s)}(x_u).$$
\end{lem}

From now on, assume $G$ is reductive. The subgroup $C_G(x_s)$ is a reductive
subgroup of $G$ of the same rank. 
Typically, it is not connected, but the identity component $C_G(x_s)^0$ contains all unipotent elements of $C_G(x_s)$.
Consequently, an element $x$ is (sub)regular in $G$ iff $x_u$ is (sub)regular in $C_G(x_s)^0$. 

For a given regular unipotent element $u$ in $G$,

\begin{lem} \label{lem:RegUniCent}
Assume that $u$ is contained in the unipotent radical 
$U$ of the Borel subgroup $B$. Then 
\begin{enumerate}
\item
$C_G(u)=Z(G)C_U(u) \isom Z(G) \times C_U(u)$,
\item
The subgroup $C_U(u)$ is connected,
\item
$C_G(u)$ is connected if and only if $Z(G)$ is connected.

\end{enumerate}
\end{lem}

We now focus on what happens in the case of conformal groups. 

\begin{lem} \label{lem:SSCentIsConnected}
Given a semisimple element $s$ in $\cgrr$, the centralizer $C_\cgrr(s)$
is a connected reductive group of the same rank as $\cgrr$.
\end{lem}

To prove this lemma, we will introduce first, a set of simple roots in $\lat_r$ and
a corresponding set of fundamental weights:

\subsection{A Set of Simple Roots} \label{ssec:SimpleRoots}
Inside $\rsys_r$, fix a set of simple roots $\a_1:=\e_1-\e_2, \dots, 
\a_{r-1}:=\e_{r-1}-\e_r,\a_r:=\e_0-\e_1-\e_2-\e_3$. Denote this set by $\rsim$ and
the corresponding set of positive roots by $\rsys_+$.
Let $\cborel$ be the corresponding  Borel subgroup of $\cgrr$ containing
$\ccartan=\lat_r\otimes\CC^*$.
From now on, we fix these three groups $\ccartan \subset \cborel \subset \cgrr$.

\subsection{Fundamental Weights}{\cite{FMdP}} \label{ssec:Lattice2}
\label{ssec:FundamentalWeights}
Corresponding to the set of simple roots, $\a_i$, for $i=1,\ldots,r$,
we introduce the fundamental weights $\overline{\omega}_i$ in $\wlat(\lat_r)$.
These will the images of elements $\omega_i\in\lat_r$:

\begin{itemize}
\item $\omega_1:=\e_0-\e_1$,
\item $\omega_2:=2\e_0-\e_1-\e_2$,
\item $\omega_i :=\e_{i+1}+\cdots+\e_r$ for $i=3,\ldots,r-1$,
\item $\omega_r:=\e_0$.
\end{itemize}

We complete both sets of elements $\alpha_i$ and $\omega_i$ to bases of
$\lat_r$ by setting $\alpha_{1+r}=e_3$ and $\omega_{r+1}=\lcan_r$. This
the unique extension with $\omega_{r+1}=\lcan_r$ and 
$\alpha_i\cdot\omega_j=\delta_{ij}$ for $i,j=1,\ldots,1+r$.

A nice consequence of the unimodularity of the lattice $\lat_r$ is 
\begin{lem} \label{lem:TSisConnected}
For any set of simple roots $S\subset \rsys_r$, $\ccartan^S:=\{t\in
\ccartan \: : \: \alpha(t)=1 \: \textnormal{ for } \: \alpha\in S \}$ is connected,
hence a torus. Its centralizer $C_\cgrr(\ccartan^S)$ is a connected reductive group
with maximal torus $\ccartan$ and connected center $\ccartan^S$.
\end{lem}
\begin{proof}
First, we will  prove that $\ccartan^S$ is connected, thus a subtorus in $\ccartan$.
The centralizer of any subtorus is a connected reductive group  \cite{HumphreysGroupBook} and  
the rest is standard. 

W.l.o.g. we  assume that $S\subset\rsim_r$.
Recall that $\ccartan\isom  \CC^* \otimes \lat_r$. Then an element
$t=\sum_{i=1}^{1+r} t_i\otimes\omega_i$ is in $\ccartan^S$ if and only if
$\alpha_j(t)=t_j=1$ for all $\alpha_j\in S$. From these equations it is
clear that $\ccartan^S\isom  \CC^* \otimes \lat_r^S$, where
$\lat_r^S:=\bigcap_{\alpha_j\in S} \ker\alpha_j$. In particular, 
$\ccartan^S$ is connected.


\end{proof}

\begin{proof}[Proof of \lemref{lem:SSCentIsConnected}]
The conformal group $\cgrr$ is a 
reductive self-dual Lie group with a connected center. Then the centralizer of any semisimple element $s$ is connected  by 
the Remark in Section 2.11 \cite{Humphreys}. W.l.o.g. assume $s\in\ccartan$. 
 Then $C_\cgrr(s)=C_\cgrr(\ccartan^S)$ where $S:=\{\alpha\in\rsim:
\alpha(s)=1 \}$. The rest now follows from the previous lemma.
\end{proof}

\begin{proof}[Proof of \thmref{thm:RegCentIsConnected}]
Given regular element $x$ with Jordan decomposition $x=x_sx_u$, 
we see that $C_{\cgrr}(x_s)$ is a connected reductive group
with connected center by the lemma above. 
Then $x_u$ is regular in the reductive group 
$C_{\cgrr}(x_s)$.
The theorem follows
by applying \lemref{lem:RegUniCent} to the pair $x_u\in C_{\cgrr}(x_s)$.
\end{proof}

\subsection{The Attached Groups} \label{ssec:GroupsAttache}
Given a del Pezzo pair $(S,Y)$ and  a marking $\varphi$, we can define more group theoretic data.
Denote the sets $\phi^{-1}(\ri(S))$ and $\phi^{-1}(\re(S))$ by
$\rsim_\varphi$ and $\R_{\varphi+}$. Let 
\begin{itemize}
\item
$\lat^\varphi:=\bigcap_{\alpha\in\rsim_\varphi} \ker\alpha$,
\item
$\ccartan^\varphi:=\ccartan^{\rsim_\varphi}=
\bigcap_{\alpha\in\rsim_\varphi} \ker\alpha$,
\item
$\cgrr_\varphi:=C_\cgrr(\ccartan^\varphi)$,
\item
$N_\varphi(\ccartan):=N_{\cgrr_\varphi}(\ccartan)$ and 
$W_\varphi:=N_\varphi/\ccartan$.
\item
Denote the  Borel subgroup of $\cgrr_\varphi$ corresponding to
the system of simple roots, $\rsim_\varphi$ by $\cborel_\varphi$ 
and its unipotent radical by $\cuni_\varphi$. 
\end{itemize}
The Levi decomposition of
$\cborel_\varphi$ is $\cborel_\varphi=\ccartan\cuni_\varphi$.
Notice that $\ccartan^\varphi=\CC^*\otimes\lat^\varphi$.
The subgroup $\ccartan^\varphi\cuni_\varphi\subset\cborel_\varphi$
is a $\cborel$-invariant (i.e. normal) subgroup of $\cborel$.

All the subgroups of $\cgrr$ above are connected. Recall that $\cgrr_\varphi$
is a reductive Lie group with maximal torus $\ccartan$ and center 
$Z(\cgrr_\varphi)=\ccartan^\varphi$. The union
$\R_{\varphi+}\sqcup(-\R_{\varphi+})$ is the root system of $\ccartan$
in $\cgrr_\varphi$, with the preferred system of simple roots, $\rsim_\varphi$.

The normalizer $N_\varphi(\ccartan)$ is $N_\cgrr(\ccartan)\cap \cgrr_\varphi$
and the Weyl subgroup $W_\varphi\subset W$ is generated by the reflections
in the roots $\in \rsim_\varphi$ (or, resp. $\in \R_{\varphi+}$)
corresponding to $-2$-curves under $\varphi$ (or resp. effective roots in
$\Pic(S)$ under $\varphi$).

When $\varphi$ is geometric, then $\cborel_\varphi=\cborel\cap\cgrr_\varphi$,
where $\cborel$ is the Borel subgroup of \secref{ssec:SimpleRoots}.

The accompanying Lie algebras may be similarly defined, for which 
all the analogous statements hold.

Let $\W{S}$ denote the group of all isometries generated by reflections in the roots of $S$. Then $\W{S}$ is the 
set all self-isometries of $\H{2}(S,\ZZ)$ preserving $[K_S]$. 

Let $\Ws{S}$ denote the subgroup of $\W{S}$ 
generated by reflections with respect to
the $(-2)$-curves in $\Pic(S)$. 
The subgroup $\Ws{S}$ can be viewed as 
the Weyl group of the singularities of the
anti-canonical model of $S$. For any marking $\varphi$,
$W_\varphi=\varphi^{-1}\Ws{S}\varphi$.




\section{Principal Bundles} \label{chap:Bundle}


We briefly discuss the relationship between vector bundles and Higgs bundles and
their principal bundle analogues. 

Constructing nontrivial vector bundles can be a hard problem depending on the base
scheme and the invariants and the  properties 
the bundles should have, such as rank, Chern classes, stability etc. In comparison, constructing Higgs bundles can be  a slightly easier
problem. 

For a given scheme $X$ and a fixed line bundle $K$ (not necessarily the canonical bundle $K_X$),
\begin{defn} 
A Higgs bundle $(V,\phi)$ on $X$ is a pair consisting
of a vector bundle $V$ and a linear map $\phi: V \ra V\otimes K$. Here, $\phi$ is 
called a $K$-valued endomorphism or  a {\em Higgs field}.  
\end{defn}

Higgs bundles were first studied by Hitchin \cite{Hitchin} 
to show that the cotangent bundle
of the  moduli space of stable principal $G$-bundles on a smooth curve $X$ of genus $>2$
is an algebraically complete integrable system. In his case,   
$K=K_X$ and the pair $(V,\phi)$ represents a point in this 
cotangent bundle. 

For Higgs bundles, the problem of construction is easier since it can be 
reformulated in terms of
the eigenvalues and the eigenspaces of $\phi$.

\subsection{Spectral Covers} \label{ssec:SpectralCovers}
In a given trivialization, one can think of $\phi$ as a square matrix. Then the characteristic equation
$\det(t\Id - \phi)=0$ defines a finite cover of $X$ in $X\times \AA^1$ contained as a divisor and the sheaf of eigenspaces on this finite cover is the sheaf $\coker(t\Id-\phi)$.
One can easily globalize these definitions and establish
a one-to-one correspondence between
Higgs Bundles $(V,\phi)$ and the spectral data $(\overline{X},\Eig)$:
\begin{equation*}
\xymatrix{
{\biggl\{
\begin{array}{c}
\textnormal{Higgs Bundles } \\
(V,\phi) \textnormal{ with } \\
\textnormal{ regular } \phi\\
\end{array}
\biggr\}} \ar@<1ex>[r]^{}
&
{\biggl\{
\begin{array}{c}
\textnormal{Spectral Data } (\overline{X},\Eig)\\
 \textnormal{ with } \overline{X}\subset \Tot(K)\\
\textnormal{ and } \Eig\in\Pic(\overline{X})\\
\end{array}
\biggr\}} \ar@<1ex>[l]_{}
}
\end{equation*}

Here
$\Tot(K)$ is the total space of the line bundle $K$, $\pi$ is the bundle projection
$\Tot(K) \ra X$  and  $\tau$ is the tautological section of $\pi^*K$ on $\Tot(K)$.  

The pair $(V,\phi)$ induces a spectral cover $\overline{X}:=\det(\Id\otimes\tau-\pi^*\phi)=0$ and the sheaf of eigenvalues $\Eig:=\coker(\Id\otimes\tau-\pi^*\phi)$. In order to control the singularities of $(\overline{X},\Eig)$, one must impose some genericity assumption on $\phi$: semisimplicity
(see \cite{BNR} for this case) or regularity. We will assume that $\phi_x$ is regular for all $x\in X$, i.e. a unique Jordan block for each eigenvalue. Notice that $\phi$
is regular everywhere if and only if  
$\Eig$ is a line bundle on $\overline{X}$. 

Conversely, starting with a divisor $\overline X \subset |n\pi^*(K)|$ s.t. it is finite and flat over $X$ and a line bundle 
$\Eig$ on $\overline X$  for some $n>0$, one can construct a rank $n$ vector bundle $V$ and a Higgs field $\phi$ by pushing the line bundle $\Eig$ and the
linear map $1\otimes \tau: \Eig \ra \Eig\otimes \pi_{|\overline X}^*(K)$ forward
by $\pi_{|\overline X}$.

These two algorithms establish a correspondence between the two sides. Of the two,
the second one is more interesting as  it enables us
to construct Higgs bundles, in particular vector bundles, from much simpler data.

\subsection{Principal Bundles} Very shortly,
a principal bundle to a vector bundle is what a group is to a vector space:
Given a scheme $X$ and a complex Lie group $G$,
\begin{defn*}
The pair $(E_G,p)$ is called a principal $G$-bundle on $X$ 
if there is a free right $G$-action on $E_G$ and $p: E_G \ra X$
is locally trivial in the etale topology so that $p(eg)=p(e)$ for all
$e\in E_G$ and $g\in G$.
\end{defn*}
The morphism $p$ is called the bundle projection. The group $G$ is called the structure group of the
bundle $E_G$. For any scheme $F$ with a left $G$-action, we can form the fiber bundle $E_G(F):=E_G\times^G F:=E_G \times F/\sim$, where $(eg,f) \sim (e,gf)$ for $e\in E_G$ and $f\in F$.
We wil denote $E_G(\gf)$ by $\Ad(E_G)$ and $E_G(G)$ by $\Aut(E_G)$ where $G$ acts on itself and
its Lie algebra $\gf$ by conjugation. 

From now on, we assume $G$ is reductive. We fix a maximal torus $T$ and a Borel subgroup $B$ 
with Levi decomposition $B=TU$. Denote the normalizer of $T$ in $G$ by $N$ and 
let $W=N/T$, the Weyl group of $T$.

\subsection{Abstract Higgs Bundles} 
In \cite{DG}, Donagi and Gaitsgory provide an analogue of the above correspondence
in the setting of principal bundles.
As expected, because of the nonlinear nature of groups, 
the steps of the solution are  much more involved. 

We will only draw a sketch of the objects involed in this new equivalence. For details, 
see \cite{DG}.
On the left side of the new equivalence is an {\em abstract} Higgs bundle:

\begin{defn} \label{defn:AbsHiggs1}
An {\em abstract Higgs bundle} is 
a pair $(E_G, \cf)$, where $E_G$ is 
a principal $G$-bundle and $\cf$ is a vector {\em subbundle} of $\Ad(E_G)$, locally the centralizer of a section $\phi$ of $\Ad(E_G)$, where $\phi_x$ is a regular element for all $x\in X$. The subbundle $\cf$
is called a {\em Higgs structure}.
\end{defn}

Notice that instead of gluing local regular sections of $\Ad(E_G)$,
the authors prefer to  glue their centralizers to give the subbundle $\cf$ of the definition. 

The right hand side is a triple consisting of again a cover $\wX$ of $X$, 
a $T$-bundle
on $\wX$ with the correct ramification behaviour, and a comparison homomorphism between $N$ and
the automorphism group of the $T$-bundle, all of which are subject to some compatibility conditions.



On the right hand side, the spectral cover $\overline X$ 
is replaced by a cameral cover $\wX$ 
and the line bundle $\Eig$ is replaced by a $T$-bundle.

\noindent {\bf Examples.} Any unramified $W$-cover is cameral. Any double cover is also cameral.
The following example is not just an example, but it is the {\em universal
example} (see Section 7 \cite{DG}).
\begin{exm}
Let $\gf=Lie(G)$ and $l=\rank(G)=\rank(\gf)$. Denote the (open) subvariety of 
regular elements (\defnref{defn:Regular}) of $\gf$ by $\gf_{reg}$.

All Cartan subalgebras are regular centralizers and conjugate under the conjugation, hence $G/N$ is the parameter space of Cartan subalgeras in $\gf$.
Because all regular centralizers are $l$-dimensional, there is a morphism sending a regular element $x\in\gf_{reg}$ to its centralizer 
$[\cf_\gf(x)\subset \gf]$ in the Grassmannian of $l$-planes in $\gf$, $\Gr(l,\gf)$.
The image is $G$-invariant under conjugation, 
contains 
a copy $G/N$, namely the subvariety of regular semisimple centralizers, and is denoted by $\GNb$.  

The variety $\GNb$ has a cameral cover $\GTb$ defined inside $\GNb\times G/B$ by$\{ (\rf, \bof) \: : \: \rf \subset \bof \} $, identifying $G/B$ with the flag variety of $G$, the parameter space of Borel subalgebras in $\gf$.

The universal Higgs bundle is the pair $(E_G, \cf)$, where $E_G$ is the trivial $G$-bundle on $\GNb$ 
and $\cf$ is the restriction of the universal $l$-plane bundle on $\Gr(l,\gf)$.
The universal cameral cover is   $\GTb \ra \GNb$ and
the universal $T$-bundle is the pullback of the tautological $T$-bundle on $G/B$ given by $G/U \ra G/B$.

 These objects are universal in the sense that any other Higgs bundle and the related cameral
 cover are obtained as their pullbacks.

\end{exm}

\begin{rem} \label{defn:AbsHiggs2}
One can show that an abstract Higgs bundle $(E_G,\cf)$ on a given scheme $X$ is
equivalent to another pair $(E_G,\sigma)$, where $E_G$ is the same principal
$G$-bundle and $\sigma: E_G \ra \overline{G/N}$ is a (right)
$G$--equivariant morphism. In fact, in \cite{DG}, the authors adopt such pairs 
$(E_G, \cf)$ as  abstract Higgs bundles.
\end{rem}

\subsection{Regularized Bundles}
In our work, we will use a version of Higgs bundle adapted to fibrations:

\begin{defn}
Given a projective flat map $p: Y \ra X$ with integral fibers,
a {\em regularized} $G$-bundle on 
$Y/X$ is a triple $(E_G,\subc_X,i)$, where
\begin{itemize}
\item $E_G$ is a principal $G$-bundle on $Y$,
\item $\subc_X$ is a vector bundle on $X$, 
and
\item The map $i: p^*\subc_X \ra \Ad(E_G)$ is a vector bundle injection, whose
image is a vector subbundle of regular centralizers, i.e. any point $x\in X$ 
has an open neighborhood $U$ (in the etale 
topology or the complex topology), so the image of $\subc_X$ under $i$ is 
the centralizer of a regular section of $\Ad(E_G)$ over $U$.  
\end{itemize}
\end{defn}

In this case, the pair $(\subc_X,i)$ will be called a {\em regularization} for
the bundle $E_G$. In other words, a regularization for the bundle $E_G$
is a Higgs structure constant along the fibers of $Y/X$. Here, we quote 
Theorem 18.5 \cite{DG}:

\begin{thmno}[Theorem 18.5 \cite{DG}]   \label{thm:FromDG}
A regularized $G$-bundle on $Y$ is the same as a triple:

\smallskip

\noindent (a) A cameral cover $\widetilde X\ra X$,

\smallskip

\noindent (b) A $W$-equivariant map $v:\widetilde X \ra \Bun_T(Y/X)$ 
(of $X$-schemes),
satisfying:

\noindent
$\alpha_i \circ v_{|D^{\alpha_i}}  =1 \in \Pic(Y/X), \forall$ simple root
$\alpha_i$, and

\smallskip

\noindent (c) An object of
$\on{Higgs}'_{\widetilde X}(X) \underset{\Tors_{T_{\widetilde X}}}\otimes \Q_v$.

\end{thmno}

Here we briefly  explain the terms of this theorem (see Section 18 \cite{DG}
for a more complete description):
$\on{Higgs}'_{\widetilde X}(X)$ encodes $T$-bundles on $\wX$ with certain ramification properties.
In general, the sheaf $T_{\widetilde X}$ is defined as a subsheaf 
of another sheaf
$\overline{T}_{\wX}$. The sheaf $\overline{T}_{\wX}$ on 
$(U\sra X)\in\Sch_{et}(X)$
is defined to be $\overline{T}_{\wX}(U)=\Hom_W({U\times_X\wX},T)$,
the group of $W$-equivariant morphisms ${U\times_X\wX}\ra T$. 
For a reductive group $G$, whose derived subgroup is simply connected, $T_{\widetilde X}=\overline{T}_{\wX}$ (p.123 \cite{DG}).

Given a sheaf of Abelian groups $\A$ on $\Sch_{et}(X)$, for $(U\sra
X)\in\Sch_{et}(X)$, denote the category of $\A|_U$-torsors on $U$ by $\Tors_\A(U)$. (Here $\Sch_{et}(X) $ denotes the big etale site over $X$, the category of all schemes over $X$ whose covering maps are surjective etale morphisms.)
The assignment $U \ra \Tors_\A(U)$ gives a sheaf of Picard categories
on $\Sch_{et}(X)$. Recall that a Picard category is a groupoid together
with a tensor category structure so that all objects are invertible.
For example, the category of $T$-bundles on the scheme $X$ is a Picard
category. Having said that, denote the sheaf defined by $U \sra \Tors_\A(U)$
by $\Tors_\A$. Letting $\A=T_{\wX}$, we get $\Tors_{T_\wX}$.

The remaining piece is $\Q_v$. For $(U\sra X)\in\Sch_{et}(X)$, $\Q_v(U)$ can
be described as the category of all possible lifts of the value map $v$
to a $T_\wY$-torsor on $Y_U:=U\times_X Y$. If $Y/X$ has a section, $\Q_v$ 
is trivial.

When $X$ is projective, the set of isomorphism classes of regularized bundle
with fixed $(\wX,v)$ is a torsor over the Abelian group $\H{1}(X,T_\wX)$ (Corollary 18.7 \cite{DG}).

We want to apply the framework of regularized bundles and a variant of this theorem to study almost regular bundles 
 on (Gorenstein) del Pezzo families.

\subsection{Almost Regular Bundles} {\em From now on, assume that all the (Gorenstein) del Pezzo surfaces are of degree $\leq 6$. 
	For families, assume the base is connected and the degree of the fibers are $\leq 6$  as well.}

Starting  with a del Pezzo pair $(S,Y)$, fix an ample divisor $H$ on $S$, a marking $\varphi$ of $S$
and fix a $\ccartan$-bundle $E_\ccartan$ a representative of $\varphi \in \Bun_\ccartan(S)$.

Via the marking, the ample divisor $H$ determines a Weyl chamber, 
hence a Borel subgroup, 
w.l.o.g. assume the Borel subgroup determined by $\varphi$  is the Borel subgroup $\cborel$ of
\secref{ssec:SimpleRoots}.

In \cite{FMdP}, the authors prove the following:

\begin{thm} \label{thm:FMdPMain}
Up to isomorphism, there is a {\em unique} $\cgrr$-bundle $E_\cgrr$ 
on  surface $S$ with the following properties:
\begin{enumerate}
\item
There is a $\cborel$-bundle $E_\cborel$ such that 
$E_\cborel/\cuni \isom E_\ccartan$ and
$E_\cgrr\isom E_\cborel \times^\cborel \cgrr$, 
where $\cuni$ is the unipotent radical of $\cborel$ and $E_\ccartan$
is a $\ccartan$-bundle with $[E_\ccartan]=\cl(\varphi)$ for a given  marking 
$\varphi$ of $S$. 
\item
The bundle $E_\cgrr$ is rigid, namely $\H{1}(S,\Ad(E_\cgrr))=0$.
\end{enumerate}

The bundle $E_\cgrr$ also has the properties:

\begin{enumerate}
\item
For every $D\in|-K_S|$, $\h{0}(D,\Ad(E_\cgrr)|D)=1+r=\rank(\cgrr)$.
\item
The bundle $E_\cgrr$ descends from $S$ to $Y$.
 \end{enumerate}
\end{thm}

\begin{defn} \label{defn:AlmostRegular1}
A $\cgrr$-bundle $E_\cgrr$ on $S$ of 
satisfying (1) and (2) of the above theorem, or its descent to $Y$, is 
 called an {\em almost regular}
$\cgrr$-bundle.   
\end{defn}

The reason why  we name these bundles almost regular is as follows:
Recall that  a semistable $G$-bundle on a smooth elliptic curve is regular if 
the infinitesimal transformations of the bundle is of dimension $\rank G$. Almost
regular bundles satisfy this property on a (Gorenstein) del Pezzo surface. Also,
for any smooth $D \in |-K_S|$ (resp. $D \in |-K_Y|$), the restriction is a regular
semistable $\cgrr$-bundle on the smooth elliptic curve $D$. However, when $D$ constains $(-2)$-curves, then the restriction to $(-2)$-curves, or equivalently the
restriction to the rational double points of $Y$ shows a different behaviour, which
will be described in detail later on (\thmref{thm:Higgs4FM}). For this reason, we have coined such bundles {\em almost} regular in analogy to the fact that surface $S$ is the blow-up of $\PP^2$ at $r$ points in {\em almost} general position.

\begin{rem}
It can be shown that the isomorphism class of almost regular $G$-bundles on $S$ is independent
of the choices of
the marking, ample divisor or the Borel subgroup chosen. 

\end{rem}

The following is a refined version of the original lemma in \cite{FMdP}:

\begin{lem}[Lemma 3.8 \cite{FMdP}] \label{lem:Lift}
Given a geometric marking $\varphi$, a corresponding $\ccartan$-bundle $E_\ccartan$ and a smooth divisor $D\in|-K_S|$, one has:
\begin{enumerate}
\item
The structure group of the bundle $e_\ccartan:=E_\ccartan|D$ reduces
to the subgroup $\ccartan^\varphi$, i.e. there is a
$\ccartan^\varphi$-bundle $e_{\ccartan^\varphi}$ s.t. 
$e_\ccartan=e_{\ccartan^\varphi} \times^{\ccartan^\varphi}\ccartan$.
The bundle $e_{\ccartan^\varphi}$ can be lifted to a $e_\creg$-bundle,
where $\creg\subset\cborel$ is a regular centralizer, whose semisimple part
is ${\ccartan^\varphi}$.
\item Let $\creg\subset\cborel_\varphi$ be any such regular centralizer and $e_\creg$
an $\creg$-bundle with the properties above. Let
$e_{\cborel_\varphi}:=\extn{e}{\creg}{\cborel_\varphi}$. Then there is a lift of the
$\ccartan$-bundle to a $\cborel_\varphi$-bundle, say $E_{\cborel_\varphi}$, so that 
$E_{\cborel_\varphi}/\cuni_\varphi\keq E_\ccartan$ and $E_{\cborel_\varphi}|D\keq\e_{\cborel_\varphi}$.
\end{enumerate}
\end{lem}


\subsection{Automorphisms of $E_\cgrr$}  
First, we touch on the endomorphisms of almost regular bundles:

\begin{lem} \label{lem:GlobalEndo}
Let $\varphi$ be a geometric marking and $E_{\cborel_\varphi}$
be a ${\cborel_\varphi}$-bundle as constructed above and
$E_\cgrr:=\extn{E}{\cborel_\varphi}{\cgrr}$. Given any subgroup 
$\cgrr_1\subset \cgrr$
containing $\cborel_\varphi$ and a $\cgrr_1$-module $\gf_2\subset \gf$ containing
$\tf^\varphi\uf_\varphi$, where $\bof_\varphi=\tf\uf_\varphi$ is the Levi decomposition for $\bof_\varphi$, set $E_{\cgrr_1}:=\extn{E}{\cgrr_\varphi}{\cgrr_1}$. 
Then the inclusion
$$\H{0}(S,E_{\cborel_\varphi}(\tf^\varphi\uf_\varphi))\subset
\H{0}(S,E_{\cgrr_1}(\gf_2))$$
is an equality.

Consequently, the following natural inclusions 
$$\H{0}(E_{\cborel_\varphi}(\tf^\varphi\uf_\varphi)) \subset
\H{0}(E_{\cborel_\varphi}(\bof_\varphi)) \subset
\H{0}(E_{\cborel_\varphi}(\gf_\varphi)) \subset
\H{0}(E_{\cborel_\varphi}(\gf))=\H{0}(\Ad E_{\cgrr})$$
are  equalities.
\end{lem}

For the definitions of the subgroups of $\cgrr$ used in the lemma,
see  \secref{ssec:GroupsAttache}.

\begin{proof}[Sketch.]
The proof is an adaptation
of Lemma 3.5 and Lemma 3.6 \cite{FMdP} to arbitrary $\cborel_\varphi$-modules
in $\gf=\Lie(\cgrr)$ combined with 
the consequences of the fact that $E_\cgrr|D$
is a regular bundle on the  smooth anticanonical divisor $D$. 
\end{proof}

\begin{rem} 
The above result indicates that the global sections of 
any of the vector bundles considered in the lemma
take their values in the Lie subalgebra $\tf^\varphi\uf_\varphi
\subset \bof_\varphi$.
\end{rem}

\label{ssec:GlobalAuto}
\begin{lem}

For all divisors $D\in|-K_S|$, the restriction map
\begin{equation} \label{eqn:RestrictionOfAuto}
\H{0}(S,\Aut(E_\cgrr)) \ra \H{0}(D,\Aut(E_\cgrr)|D)
\end{equation}
is  injective. 
When $D$ is smooth, it is  an isomorphism.
The group of global automorphisms, $\H{0}(S,\Aut(E_\cgrr))$, of an almost regular bundle 
$E_\cgrr$ is isomorphic to a regular centralizer in $\cgrr$ and is, therefore, connected.
\end{lem}


\begin{proof}

Recall that  the restriction map 
\begin{equation} \label{eqn:RestrictionOfEndo}
\H{0}(S,\Ad(E_\cgrr)) \ra \H{0}(D,\Ad(E_\cgrr)|D)
\end{equation} 
is an isomorphism
for all $D\in |-K_S|$. 
Because $\H{0}(\cdot,\Ad(E_\cgrr))$ is the tangent Lie algebra of the Lie group
$\H{0}(\cdot,\Aut(E_\cgrr))$, we conclude that \eqnref{eqn:RestrictionOfAuto} is
an isomorphism of identity components.

We now show that it is injective for arbitrary $D$:
Pick $g\in\H{0}(S,\Aut(E_\cgrr))$ so that $g|D=\Id$ on $D$. By differentiating $g$,
we see that $dg|D=0$, i.e. $dg\in\H{0}(S,\Ad(E_\cgrr)(-D))$. This group is trivial, hence $dg$ vanishes identically on $S$. From this it follows that  $g$ is constant. 
Since
$g|D=\Id$ on $D$, we see that $g=\Id$ on $S$. Hence, the restriction map
on the global automorphisms of $E_\cgrr$
is injective.

When $D$ is smooth, it is an isomorphism:
$E_\cgrr|D$ is a regular bundle on
the elliptic curve $D$ and $\H{0}(D,\Aut(E_\cgrr)|D)$ is isomorphic to a regular centralizer in $\cgrr$ (\cite{FM1}), which is connected ( \thmref{thm:RegCentIsConnected}). This proves that \eqnref{eqn:RestrictionOfAuto} is
an isomorphism for smooth $D$. We conclude that $\H{0}(S,\Ad(E_\cgrr))$ is isomorphic to
a regular centralizer in $\cgrr$ in $\cgrr$ and connected as well.

\end{proof}

\subsection{Families} We define families of del Pezzo surfaces, almost regular bundles on them and prove a lemma comparing the cohomology line bundles of a
family of del Pezzo surfaces and anticanonical divisors:
 
\begin{defn}
A family of (Gorenstein) del Pezzo surfaces is a flat projective map 
$p: (\YY) \SS \ra \XX$ whose geometric fibers are (Gorenstein) del Pezzo surfaces.
\end{defn}

\begin{defn} \label{defn:FMonFamily} \label{defn:AlmostRegular2}
An almost regular bundle $E_\cgrr$ on a family of (Gorenstein) del Pezzo surfaces 
is
a principal bundle $E_\cgrr$ on the total space  whose restriction to each (geometric) fiber is almost regular.
\end{defn}

\begin{rem} Recall that on a single surface, all such bundles are isomorphic.
In fact, we will see that any two such bundles are also isomorphic locally with respect to the base.
\end{rem}

Given a family of del Pezzo surfaces , $p: S \ra X$, of degree $9-r$, 
a map $\varphi: \lat_S \ra \Pic(S/X)$, is a relative marking if the restriction of $\varphi$ to each geometric
fiber is a marking. Here, $\lat_S$ is the locally constant sheaf on $S$ attached to $\lat=\lat_r$.

\begin{lem} \label{lem:CompareCohomology}
For a given root $\alpha$, let $L_\alpha$ be a line bundle corresponding to $\alpha$ and $D$ a relative anticanonical section with integral fibers. Then

\begin{enumerate}
\item $R^ip_*L_\alpha=R^i{p_{D}}_*L_\alpha=0$ for $i\geq 2$.
\item For roots $\alpha$ positive w.r.t. $\varphi$, the restriction maps $R^ip_*L_\alpha \sra R^i{p_{D}}_*L_\alpha$ are isomorphisms for $i=0,1$ and
\item The restriction maps $\H{j}(S,L_\alpha) \sra \H{j}(D,L_\alpha|D)$ are isomorphisms for $j\geq 0$. 

\end{enumerate}

\end{lem}

\begin{proof}

1) and 2) follow from Lemma 3.2 and 3.3 of \cite{FMdP}.
To prove 3), we see that by 1), the Leray spectral sequence for both $L_\alpha$ on $S$ 
and $L_\alpha|D$ on $D$ degenerates at the $E_2$-term, resulting in a long exact sequence for each.
The restriction map induces a commutative diagram whose rows are  these exact sequences.
Applying 2) and the five lemma to this commutative diagram yields 3).
\end{proof}

\subsection{Local Models}
\label{ssec:LocalModel}
The following is a consequence of the previous lemma 
(in \secref{ssec:GlobalAuto}):
\begin{thm} \label{thm:IsomLocalModel}
Let $p:\YY \ra \XX$ be a family of  (Gorenstein) del Pezzo surfaces. 
Suppose ${E_\cgrr}$ is an almost regular bundle on $\YY$.
Then the sheaves $p_*\Ad(E_\cgrr)$ and  $p_*\Aut(E_\cgrr)$ are locally trivial
on $\XX$. 

Consequently, if $E_\cgrr^i$, $i=1,2$ are two such bundles, then 
$p_*\Isom(E_\cgrr^1,E_\cgrr^2)$ is locally trivial on $\XX$. 
Therefore, given a point $x\in \XX$, there is an open (etale) neighborhood
$\UU$ of $x$, for 
which the restriction of the two bundles $E_\cgrr^i$ to $\UU \times_\XX \YY$
are isomorphic to each other.
\end{thm}

\begin{proof}
The sheaf $p_*\Ad(E_\cgrr)$ is locally free since for any geometric
fibre $F$,
\\$h^{0}(F, \Ad(E_\cgrr)|F) = 1+r$.  Combining this with the fact that the group \\$\H{0}(F,\Aut(E_\cgrr)|F)$ 
is connected, we see that $p_*\Aut(E_G)$
is a locally trivial group scheme with $(1+r)$--dimensional connected fibers.

Given two 
almost regular bundles $E^i_\cgrr$, $i=1,2$, the sheaf 
$p_*\Isom(E^1_\cgrr,E^2_\cgrr)$ is
a torsor over $p_*\Aut(E^1_\cgrr)$, hence locally trivial. The final
conclusion follows from this fact.
\end{proof}

We make the following definition:
\begin{defn}
We will say that two principal $G$-bundles on a fibration $Y/X$ are locally isomorphic with
respect to the base $X$, if for a given point in $X$, there is 
a neighborhood $U$ of the point so that the pullback of the two bundles 
on $U \times_X Y$ are isomorphic.
\end{defn}

The above theorem shows that two almost regular bundles on a family
of del Pezzo surfaces are locally isomorphic w.r.to the base.
However, notice
that the neighborhood $\UU$ is not universal in any sense; it depends
on the two bundles in consideration.



\subsection{The Higgs Structure of Almost Regular Bundles}



\begin{lem} \label{lem:SmoothAntiKDiv}
Given a del Pezzo surface $S$ and a point $s$ on $S$ not contained in any $(-2)$-curve, there exists a smooth divisor $D\in |-K_S|$ containing $s$.
\end{lem}
\begin{proof}
By Bertini's Theorem.
\end{proof}

\subsection{Regularity on Elliptic Curves}
 Given a reductive Lie group, G, recall that
\begin{defn} \label{defn:RegularElliptic}
A principal $G$-bundle $E_G$ on a smooth curve $D$ is {\em semistable}
if its adjoint bundle $\Ad(E_G)$ is semistable. It  is {\em
regular} if $\h{0}(D,\Ad(E_G))=\rank G$.

\end{defn}

\begin{lem}
Let $V$ be a semistable vector bundle of degree $0$ on a smooth curve
$D$. Then the evaluation map 
$$\H{0}(D,V)\otimes \O_D \ra V$$
is an injective map of vector bundles.
\end{lem}
\begin{proof}
It is enough to show that any nonzero section $s$ of $V$ is has no zeroes.

If $s$ vanishes along a divisor $E$ on $D$, we see that the line bundle 
$\O(E)$ is a line sub bundle of $V$. Since $V$ is semistable of slope
$0$, we conclude that there is no such $E$. In other words, $s$ has no zeroes, proving that the evaluation map is injective. 
\end{proof}

\begin{cor}
Let $D$ be a smooth elliptic curve and $E_G$ a regular bundle on $D$.
Then
$$\H{0}(D,\Ad(E_G))\otimes  \O_D \ra \Ad(E_G)$$
is an injection of vector bundles.
\end{cor}

\subsection{Regularity on the Surface $S$}

\newcommand{\ev}{{\textnormal{ev}}}
\begin{thm} \label{thm:Higgs4FM}
An almost regular $\cgrr$-bundle $E_\cgrr$ on $S$ determines a Higgs structure
on $S^0$, namely 
$\ev^0: \H{0}(S,\Ad(E_\cgrr))\otimes \O_{S^0} \ra
\Ad(E_\cgrr)|S^0$. In other words,
\begin{enumerate}
\item
The evaluation map $\ev^0$ is of maximal rank ($=1+r$) on $S^0$,
\item
The image of $\ev^0$ provides a Higgs structure for $E_\cgrr$.

\noindent
Furthermore,
\item

Although, the evaluation map  $\ev: \H{0}(S,\Ad(E_\cgrr)\otimes \O_{S} \ra
\Ad(E_\cgrr)$ is well--defined, the image of $\ev$ along $-2$--curves, has no
regular elements in it. In other words,
the Higgs structure given by $\H{0}(S,\Ad(E_\cgrr))\otimes \O_{S^0}$ and $\ev^0$
 does not extend any further.
\end{enumerate}
\end{thm}
\begin{proof}
To prove (1), recall that given a point $s\in S^0$, there is 
a smooth divisor $D\in|-K_S|$  by
 \lemref{lem:SmoothAntiKDiv}.  By the corollary above, 
it follows that the
evaluation map is of maximal rank along $D$, hence $\ev$ is of maximal 
 for every $s\in S^0$.

To see (2) and (3),
 fix a geometric marking $\varphi$. Notice that 
$\tf^\varphi$ is a trivial $\cborel_\varphi$-module,
$\tf^\varphi\uf_\varphi\subset \bof_\varphi$ 
is a $\cborel_\varphi$-module and it splits
since
$\tf^\varphi\subset \zf(\bof_\varphi)$. 
 It follows that 
$E_{\cborel_\varphi }(\tf^\varphi)$ is a trivial vector bundle and 
the nontrivial part
of the bundle $E_{\cborel_\varphi }(\tf^\varphi\uf_\varphi)$ is
$E_{\cborel_\varphi }(\uf_\varphi)$. By \lemref{lem:GlobalEndo}, 
we have 
$$
\begin{array}{rcl}
\H{0}(S,\E_\cgrr(\gf)) & \isom &
\H{0}(S,E_{\cborel_\varphi}(\tf^\varphi\oplus\uf_\varphi)) \\
 & \isom &
\H{0}(S,E_{\cborel_\varphi}(\tf^\varphi))
\oplus
\H{0}(S,E_{\cborel_\varphi}(\uf_\varphi)). \\
\end{array}$$


Any section $x\in\H{0}(S,E_{\cborel_\varphi}(\tf^\varphi\oplus\uf_\varphi))$
has a Jordan decomposition with $x=x_s+x_n$ s.t.
 $x_s\in\H{0}(S,E_{\cborel_\varphi}(\tf^\varphi))$
 and
 $x_n\in\H{0}(S,E_{\cborel_\varphi}(\uf_\varphi))$.
 A section $x$ will be  regular at a point $p\in S$ iff
$\cf_\gf(x_s(p))=\gf_\varphi$ and $x_n(p)$ is a regular nilpotent in 
$\cf_\gf(x_s(p))$.
Since $x_s$ is constant, the first condition holds at one point $p\in S$ 
if and only 
if it holds throughout  $S$. From now on, assume that this is the case.

To check the regularity of the nilpotent section $x_n$ at a point $p$,
it is necessary and sufficient to see if all the $\alpha$-component of $x_n(p)$
are nonzero for 
simple roots $\alpha\in\Sigma_\varphi:=\varphi^{-1}(\ri(S))$:
Let $\uf_\varphi^2:=[\uf_\varphi,\uf_\varphi]$, then there is a
$\cborel_\varphi$-module sequence:
$$ 0 \ra \uf_\varphi^2 \ra \uf_\varphi \ra
\oplus_{\alpha\in\Sigma_\varphi}\gf_\alpha \ra 0.$$
Denote the projection map $\uf_\varphi \ra \gf_\alpha$ by $\pi_\alpha$.
A nilpotent valued section $x_n$ at a point $p$ is regular iff  
$\pi_\alpha(x_n(p))\neq 0$ for all simple roots $\alpha\in\Sigma_\varphi$.

Fix a smooth divisor $D\in|-K_S|$ and a section $x=x_s+x_n$ such that
$x|D$ is regular. This is possible since $E_\cgrr|D$ is regular and
$\H{0}(D,\Ad(E_\cgrr)|D)$ is a regular centralizer contained in $\tf^\varphi\uf_\varphi
\subset\gf$.

For $\alpha\in\Sigma_\varphi$, the roots  corresponding to $-2$--curves under
$\varphi$, we have $E_\cborel(\gf_\alpha)|D\isom \O_D$. It follows that
$\pi_\alpha(x_n)|D$
is constant and $\neq 0$. 
Since $\H{0}(S,E_\cborel(\gf_\alpha)) \isom 
\H{0}(D,E_\cborel(\gf_\alpha)|D) \isom \CC$, it follows that $\pi_\alpha(x_n)$
is a nontrivial section  of the line bundle $E_\cborel(\gf_\alpha)$. Thus
$\div(\pi_\alpha(x_n))$ is exactly the $-2$--curve $C_\alpha$, where
$E_\cborel(\gf_\alpha)=\O_S(C_\alpha)$. Away from $C_\alpha$,
the section $\pi_\alpha(x_n)\in E_\cborel(\gf_\alpha)$ has no zeroes,
which proves that $x_n$ is a regular in $\gf_\varphi$ and $x$ is regular in $\gf$ on $S^0$. By the same argument, we see that $\alpha$-component of the section $x$
vanishes along the $(-2)$-curve $C_\alpha$.
The former proves that the image of $\ev^0$ is the centralizer of such a section $x$ and hence is a regular centralizer subbundle, i.e. a Higgs structure for $E_\cgrr$
proving (2). The latter proves (3) that this Higgs structure cannot be extended to 
the union of $(-2)$-curves.


\end{proof}

We adopt a more relaxed definition of regularization:
\begin{defn}
Given a projective flat map $p: Y \ra X$ with integral fibers and open subscheme $Y^0 \subset Y$ with
$\codim(Y - Y^0)\geq 2$,
a {\em regularized} $G$-bundle on 
$Y^0$ is a triple $(E_G,\subc_X,i)$, where
\begin{itemize}
\item $E_G$ is a principal $G$-bundle on $Y^0$,
\item $\subc_X$ is a vector bundle on $X$, 
and
\item The map $i: p^*\subc_X \ra \Ad(E_G)$ is an injection of coherent sheaves,
which is of maximal rank over the open subscheme $Y^0$ whose
image is a vector subbundle of regular centralizers.

\end{itemize}
\end{defn}

\begin{thm} \label{thm:UniqueReg}
An almost regular bundle $E_\cgrr$ on $p:\YY \ra \XX$ admits a unique regularization,
namely $p_*\Ad(E_\cgrr)$.
\end{thm}
\begin{proof}
Keeping the notation of the above definition,
observe that $i: p^*\subc_X \sra
\Ad(E_G)$ factors through $p^*p_*\Ad(E_G)$, i.e.  
$p^*\subc_\XX \ra p^*p_*\Ad(E_G) \ra \Ad(E_G).$
If we start with an almost regular bundle $E_\cgrr$ on a family $p: \YY \ra \XX$, then 
$E_\cgrr$ is regularized by $p_*\Ad(E_\cgrr)$ (this 
follows from \thmref{thm:Higgs4FM}). The previous observation shows that this is the unique
regularization possible.
\end{proof}





 \newcommand{\wE}{{\overline E}}
 
 From now on, assume the Gorenstein del Pezzo family $\YY/\XX$ has a section.
 Then the sheaf $Q_v$ is trivial.
 
 The following is our main theorem.
 
 \begin{thm} \label{thm:MainGeneral}
  The following datum 
\begin{enumerate}
\item A cameral resolution
$$\xymatrix{
  \SS \ar[r]^r \ar[d]_s 
&  \wYY\ar[d]_{\tilde{p}} \ar[r]^{\widetilde c} & \YY \ar[d]_p \\
 \wXX \ar@{=}[r] & \wXX \ar[r]^c & \XX\\ 
}$$
with normal cameral cover $\wXX$,
\item A marking $\Phi:\lat_\SS \ra \Pic(\SS/\wXX)$, 
\item An object from $\higgsp{X}$ 

\end{enumerate}
determines  a unique almost regular $\cgrr$-bundle
on $\YY$.

\end{thm}
 
When there is a cameral resolution as above, there exists such markings. By Proposition 6.8.2 \cite{EGA-IV},
it follows that $\wYY$ (resp. $\SS$) since normal since $\wYY$ has normal  (resp. $\SS$ has smooth) fibers over the normal scheme $\wXX$. Since $\SS \ra \wYY$
induces an isomorphism between the $\wYY^0$ and its preimage, we identify these two.

\begin{rem} In the proof of the above statement, we will use 
Theorem 18.5 \cite{DG} (cited as \thmref{thm:FromDG}). However, the careful
reader will notice that we cannot possibly apply this theorem to our setup as it is
unless we impose that $\YY\ra \XX$ is a smooth morphism. Therefore, without
ever stating a precise version of this theorem applicable directly in our case, 
we will informally discuss how it is possible to apply this theorem under the assumptions  we presented  above. In short,
this is possible by amending the reasoning for Theorem 18.5 
by requiring that $\wXX$ is normal and subsequently, invoking Hartogs' Lemma on $\SS$. 
 
For this theorem, one needs a projective flat morphism with  integral fibers.
The assumptions, projectivity and integral fibers are used for two purposes.
The first is the existence of the relative Picard scheme. One can get around this 
by dealing with the global sections of the relative Picard sheaf. The second is to show
that sections, automorphisms etc. of some $T$- (and related $\CC^*$-) are bundles are constant (possibly along the fibers) using the properness. Such arguments used in the proof of \thmref{thm:FromDG}  fails for $\YY^0$. However, we can reach
the same conclusions by invoking Hartogs' Lemma for
such objects on $\SS$:
For instance, as a part of the proof of \thmref{thm:FromDG}
one has to show that the certain automorphism of $E_\ccartan^0$ is trivial.
Starting with the automorphisms on $\YY^0$, we see by Hartogs' Lemma that they
extend to automorphisms of $E_\ccartan$ on $\SS$. (Notice that we only extend 
sections, not the objects. Therefore the existence of $E_\ccartan$ a priori to the construction is crucial.) However, $\SS$ over $\XX$ is projective with integral
fibers and hence the rest of the arguments runs as before. Notice that  the resulting 
bundles are regularized $\cgrr$-bundles on $\YY^0$, rather than on $\YY$.
\end{rem}
 
 \begin{proof}
 We prove this result in several steps. 

The first step is to build a regularized $\cgrr$-bundle $E_\cgrr^0$ on $\YY^0$ by slightly altering 
the arguments of Section 18 \cite{DG} used to prove Theorem 18.5\cite{DG} (or, \thmref{thm:FromDG}). The second step is to show that the pullback of this bundle to $\wYY^0$ extends to an almost
regular bundle $\wE_\cgrr$ on $\SS$.  The final step is to show that this bundle descends
from $\SS$ to $\YY$.

We do the final step first. Assume that there is an almost regular $\cgrr$-bundle
$\wE_\cgrr$ on $\SS$. It descends to $\wYY$ by Lemma 3.4 \cite{FMdP}. Denote this bundle by $\widetilde{E}_\cgrr$ on $\wYY$.
The second descent is from $\wYY$ to $\YY$: Pullback by $w\in W$ induces 
an automorphism of $\widetilde{E}_\cgrr^0$. Applying Hartogs' Lemma to automorphisms of  $\widetilde{E}_\cgrr$, we see that this induced automorphism 
extends to $\wYY$. It is clear that such automorphisms satisfy the descent condition for $\widetilde{E}_\cgrr$ by the density of $\wYY^0$ in $\wYY$. Consequently, $\widetilde{E}_\cgrr$ descends. 

We proceed with the first step. 
First, pick a $\ccartan$-bundle $E_\ccartan$
representing $\Phi$ on $\SS$. Restricting $E_\ccartan$ to $\wYY^0$ determines 
a value map $\Phi^0$. Now the three piece datum $(\wXX, \Phi^0, \textnormal{Object} \in \higgsp{X})$ determine a regularized $\cgrr$-bundle on $\wYY$ by the above remark.

\newcommand{\we}{{\widetilde e}}
In the second step, we will extend $\wE_\cgrr^0$ to $\SS$. This will be done
in an inductive manner following \cite{FMdP}. Notice that
for any etale open $\UU  \ra \XX$, the extension of the bundle $E_\cgrr^0$ 
from $\wYY^0 \times_\UU\XX$ to $\SS \times_\UU \XX$ is unique. Therefore,
it is sufficient to extend the bundle $E_\cgrr^0$ locally with respect to $\XX$.
Let $\DD \in |-K_{\YY/\XX}|$ have
integral fibers. Such a $\DD$ exists after possibly shrinking $\XX$.

The restriction of $(\wXX, \Phi^0, \textnormal{Object} \in \higgsp{X} \otimes \Q_{\Phi^0})$ to $\DD$ produces a
regularized $\cgrr$-bundle on $\DD$. (In fact, the restriction to a given smooth fiber is regular.) Recall that on the cameral cover $\wDD$, there is a $\cborel$-bundle,
$\we_\cborel$, s.t. $\we_\cgrr=\extn{\we}{\cborel}{\cgrr}$, 
$\we_\cborel/\cuni\cong E_\ccartan|\wDD$, where $\we_\cgrr$ is the pullback
of $e_\cgrr$ to $\wDD$. 
Also, there exists a $\cborel$-bundle $\wE_\cborel^0$ on $\wYY^0$ with
 $\wE_\cgrr^0=\extn{{\wE^0}}{\cborel}{\cgrr}$ and $\wE^0_\cborel/\cuni\isom E_\ccartan$ on $\YY^0$.

We will extend $\wE^0_\cborel$ to $\SS$. This will be done as follows: We will
lift the $\ccartan$-bundle $E_\ccartan$ to a $\cborel$-bundle and will prove
that the two $\cborel$-bundle (resp. the associated $\cgrr$-bundles) are isomophic
on $\wYY^0$.  This proves that $\wE^0_\cborel$ (resp. $\wE^0_\cgrr)$ extends to $\SS$.

The lifting problem from $\ccartan$ to $\cborel$ will be carried out in an inductive manner, where for
each step, one has  a lifting problem for $0\sra N \sra G \sra H \sra 0$ with $N$ normal and Abelian.
Two fundamental questions are: Can one a lift a $H$-bundle on $X$ 
to a $G$-bundle ?
If so, how many lifts are there ? The group $H$ acts
on $N$. Given $E_H$ one can form $E_H(N)$ and the  sequence 
$$\xymatrix{H^1(X, E_H(N)) \ar[r]  &  
H^1(X, G) \ar[r]  &
 H^1(X, H) \ar[r]^{\delta}  &            
H^2(X, E_H(N )) }.$$
The bundle $E_H$ has a lift to a $G$-bundle iff the obstruction $\delta[E_H] \in H^2(X, E_H(N ))$ vanishes.
In which case, $H^1(X, E_H(N))$ acts transitively on the set of all such lifts. For details, see \cite{ManinGauge}.

Fix a decreasing filtration
of the unipotent group $\UU$, $\cuni=\cuni_0 \supset 
\cuni_1 \supset \cdots
\supset \cuni_N \supset \cuni_{N+1}=\{ 1 \}$ by normal $\ccartan$-invariant 
subgroups
such that the quotient $\cuni_j/\cuni_{j+1}$ is contained in the center of 
$\cuni/
\cuni_{j+1}$ for all $j$. Then the quotient group $\cuni_j/\cuni_{j+1}$ is a 
direct sum of root spaces $\gf_\alpha$. Denote the index set of the roots for
the $\cuni_j/\cuni_{j+1}$
by $\R_j \subset \R$, then set
$\cuni^j:=\cuni_j/\cuni_{j+1}=\bigoplus_{\alpha\in\R_j} \gf_\alpha.$ Also set $\cborel_j:=\cborel/\cuni_j$
and define $\we_{\cborel_j}$ to be $\we_\cborel/\cuni_j$. 
The bundles $\we_{\cborel_j}$ 
will be used as a framing along the genus $1$ fibration $\wDD/\wXX$.

We do induction on $j$. 
As the initial step, set $\wE_{\cborel_0}:=E_\ccartan$.
The induction step is to lift a
${\cborel_j}$--bundle
$\wE_{\cborel_j}$ s.t. $\wE_{\cborel_j}|\wDD \isom
\we_{\cborel_j}$ to a $\cborel_{j+1}$--bundle 
$\wE_{\cborel_{j+1} } $ s.t. 
$\wE_{ \cborel_{j+1} }|\wDD\isom \we_{\cborel_{j+1}}.$

For $\wE_{\cborel_j}$ and $\we_{\cborel_j}$, $\cborel_j$-bundles on $\SS$ and $\wDD$ respectively, we have the following diagram for the lifting problem

\begin{equation}
\xymatrix{
\H{1}(\SS, \wE_{\cborel_j}(\cuni_j) ) \ar[r] \ar[d]   &  
\H{1}(\SS, \cborel_{j+1}) \ar[r] \ar[d]   &  
\H{1}(\SS, \cborel_{j}) \ar[r] \ar[d] &  
\H{2}(\SS, \wE_{\cborel_j}(\cuni_j) )  \ar[d]  \\
\H{1}(\wDD, \we_{\cborel_j}(\cuni_j) ) \ar[r]   &  
\H{1}(\wDD, \cborel_{j+1}) \ar[r]  &  
\H{1}(\wDD, \cborel_{j}) \ar[r]  &  
\H{2}(\wDD, \we_{\cborel_j}(\cuni_j) )   \\
}
\end{equation}
where the downward arrows are the restriction maps from $\SS$ to $\wDD$. 
 
The vector bundle  $\wE_{\cborel_j}(\cuni_j)$ is direct sum of line bundles, $\wE_\ccartan(\gf_\alpha)$, corresponding to 
the {\em positive} roots $\alpha \in \R_j$. By \lemref{lem:CompareCohomology}, it follows that
the first and the last vertical arrows are isomorphisms. 

Recall that the obstruction to lift $\wE_{\cborel_j}$ to a $\cborel_{j+1}$-bundle restricts to the obstruction
to the lift $\we_{\cborel_j}$. The obstruction on $\wDD$  vanishes since there is already a lift,
namely $\we_{\cborel_{j+1}}$ as prescribed above. Therefore we see that  $\wE_{\cborel_j}$ also lifts.
The first cohomology acts transitively on the set of all such lifts.
Since the first cohomology groups are isomorphic for $\SS$ and $\wDD$, we can
choose the lift $\wE_{\cborel_{j+1}}$ so that
$\wE_{\cborel_{j+1}}|\wDD \isom \we_{\cborel_{j+1}}$.
This concludes the induction process.

Once this is done, set $\wE_\cgrr:=\wE_\cborel\times^\cborel \cgrr$. 
Then by construction, $\wE_\cgrr$ is isomorphic to $\widetilde{E}_\cgrr^0$ on $\wYY^0$.
Since $\widetilde{E}_\cgrr^0$ is isomorphic to a bundle which extends to $\SS$, we
conclude that it  extends to $\SS$.
 \end{proof}
 
 \newcommand{\ASheaf}{{\mathscr A}}
Using the above theorem, we conclude that on a given (Gorenstein) del Pezzo fibration $\YY/\XX$ with normal $\wXX$, one can 
construct an Abelian sheaf, $\ASheaf$, whose cohomology groups 
 $\H{2}(\XX,\ASheaf)$ and  $\H{1}(\XX,\ASheaf)$ capture the basic nature of
 almost regular bundles on $\YY$: The obstruction to the existence of an almost
 regular bundle lives in $\H{2}(\XX,\ASheaf)$. There is such a bundle if and only if
 this obstruction vanishes. Moreover, $\H{1}(\XX,\ASheaf)$
acts transitively on the isomorphism classes of almost regular bundles and simply transitively if
$\XX$ is projective. Here, $\ASheaf$ is the sheaf of automorphisms of almost regular bundles along the fibers of $\YY/\XX$. It is a sheaf of abelian groups on $\XX$. Its existence is enabled by the above theorem, the fact that almost
regular bundles are locally isomorphic w.r.t. the base (\thmref{thm:IsomLocalModel}) and the fact that the automorphism groups
along the fibers are abelian groups (\secref{ssec:GlobalAuto}).

We see that once the cameral cover $\wXX$ is normal, and there is a simultaneous resolution $\SS$ above $\wXX$, then the bundles corresponding to the pair $(\wXX, \Phi)$ can all be extended to 
$\YY$ and the isomorphism classes of these almost regular bundles are a homogeneous space
under the $\H{1}(\XX,\ASheaf)$ action.

This is in aggrement with the results of \cite{DG}:
The set of isomorphism classes of regularized $\cgrr$-bundles on $\YY^0$ corresponding
to the pair $(\wXX, \Phi^0)$ 
is a homogenous space under the $\H{1}(\XX,\ccartan_\wXX)$ action.
Moreover, the two abelian sheaves $\ccartan_\wXX$
and $\ASheaf$ are canonically isomorphic by Theorem 11.6 \cite{DG}.

By \thmref{thm:CurveCamCover} and the above discussion, we conclude that
\begin{thm} \label{thm:MainCurve}
Given a Gorenstein del Pezzo family $\YY/\XX$ whose base $\XX$ is a smooth projective curve 
and total space $\YY$ a smooth threefold, then there 
exists
a smooth cameral cover $\wXX$ together with a simultaneous resolution
$\SS$. For a fixed relative marking $\Phi$,
the isomorphism classes of almost regular bundles  are a torsor under $\H{1}(\XX,\ASheaf)$. 
\end{thm}
\noindent {\bf Note.} In this case, $\H{2}(\XX,\ASheaf)=0$ by \cite{FM1}.

In Proposition 5.5 \cite{FM1}, the authors observe that $\H{1}(\XX,\ASheaf)$ can be related to the intermediate Jacobian
of the threefold in the following way:
Let $J^3(\YY)$ be the intermediate Jacobian of $\YY$ and $\DD$ be a smooth (relative) anticanonical divisor in $\YY$. Then define the relative intermediate Jacobian $J^3(\YY/\XX)$ as the kernel
of the morphism $J^3(\YY) \ra J(\XX)$ induced from $\H{*}(\YY) \ra \H{*}(\DD) \ra \H{*-2}(\XX)$.
Also set
$$\H{2,2}_0(\YY,\ZZ):=\ker( \: \H{4}(\YY, \ZZ) \ra \H{2}(\XX,\ZZ) \: ) / \ZZ [\YY_x]. $$
where $\YY_x$ is a general fiber. Then $\H{2,2}_0(\YY,\ZZ)$ is a finite group in general and
$\H{1}(\XX,\ASheaf)$ can be fitted into an exact sequence:

$$ 0 \ra J^3(\YY/\XX) \ra \H{1}(\XX,\ASheaf) \ra \H{2,2}_0(\YY,\ZZ) \ra 0.$$

This result generalizes Kanev's earlier results for fibrations for which
 $\XX=\PP^1$, $4\leq r \leq 7$, and the fibers have at worst a single 
ordinary node ($A_1$-singularity).

\bibliographystyle{plain}
\bibliography{almostregular}

\end{document}